\documentclass[reqno,12pt,letterpaper]{amsart}
\usepackage[proof]{macros}

\usepackage{enumerate}
\usepackage{subfig,caption}
\usepackage{graphicx}				
\usepackage{amssymb}
\usepackage{amsmath}
\usepackage{amsthm,amsfonts,bm,bbm}
\usepackage{mathtools}
\mathtoolsset{showonlyrefs=true}
\usepackage{mathrsfs}
\usepackage{xcolor}
\usepackage{enumerate,enumitem}
\usepackage{leftidx}
\counterwithout{equation}{section}
\usepackage{marginnote}
\title[Quantitative unique continuation]{Quantitative unique continuation on conic manifolds}

\author{Ruoyu P. T. Wang}
\address{Department of Mathematics, University College London, London, WC1H 0AY, United Kingdom}
\email[Ruoyu~P.~T.~Wang]{ruoyu.wang@ucl.ac.uk}

\date{}

\renewcommand{\supp}{\text{supp }}

\newcommand{\abs}[1]{\left|#1\right|}

\newcommand{\bigo}{\mathcal{O}}

\begin{document}

\begin{abstract}
This expository note, written for the proceedings of ICCM 2023, presents recent work \cite{wan20}. We particularly prove an Carleman estimate on conic manifolds, using a multiple-weight Carleman argument.
\end{abstract}

\maketitle

\section{Introduction}
Let $M$ be a smooth compact manifold without boundary, and $\Delta$ be the nonnegative Laplace--Beltrami operator on $M$. We consider the unique continuation problem for the eigenfunction of $\Delta$:\let\thefootnote\relax\footnote{In memory of Steve Zelditch.}
\begin{question}[Unique continuation property]
Fix $\Omega\subset M$ open. Is it possible for any nontrivial eigenfunction $u_h\in C^\infty(M)$ to vanish on $\Omega$, that is, 
\begin{equation}
(h^2\Delta-1)u_h=0, \ u_h\neq 0, \ u_h|_{\Omega}=0?
\end{equation}

\end{question}
When $M$ is real-analytic, the Holmgren uniqueness principle implies $u_h$ is analytic and thus cannot vanish on any open sets. When $M$ is merely smooth, it is not easy to deduce such strong properties from that $u_h$ is only smooth. The unique continuation problem was answered by Carleman in \cite{car39}:
\begin{proposition}[Carleman estimates]
Let $M$ be compact and smooth, and $\Omega\subset M$ be an open set. Then any eigenfunction of $\Delta$ satisfies the unique continuation property:
\begin{equation}
(h^2\Delta-1)u_h =0, \ u_h|_{\Omega}=0 \Rightarrow u_h\equiv 0. 
\end{equation}
Moreover, we have the quantitative estimate, called the \emph{Carleman estimate}, 
\begin{equation}\label{e1}
\|u\|_{L^2(M)}\le e^{C/h}\left(\|(h^2\Delta-1)u\|_{L^2(M)}+\|u\|_{L^2(\Omega)}\right)
\end{equation}
uniformly for $u\in H^2(M)$, where $C=C(M, \Omega)$ does not depend on $h$. 
\end{proposition}
The Carleman estimate answered the question of unique continuation: eigenfunctions $u_h$ not only cannot vanish in open sets, but also cannot be exponentially small in open sets: we have
\begin{equation}
\|u_h\|_{L^2(\Omega)}\ge e^{-C/h}\|u_h\|_{L^2(M)},
\end{equation}
for $(h^2\Delta-1)u_h=0$. See further results on Carleman estimates in \cite{leb93,lr95,lr97,bur98,ler19}.

The goal of the manuscript \cite{wan20} is generalise the Carleman estimates onto noncompact manifolds. One first observes that one must impose stronger assumptions over the set $\Omega$: estimate \eqref{e1} cannot hold for any precompact open set $\Omega$. Indeed, let $M=\mathbb{R}^d$ and $u_h=\exp(ih^{-1}\xi.x)$ be the plane waves for unit covector $\xi$. Then $(h^2\Delta-1)u_h=0$, $\|u\|_{L^2(\Omega)}=\abs{\Omega}$ equal to the finite volume of $\Omega$, but $\|u\|_{L^2(\mathbb{R}^d)}$ is unbounded. 

A natural assumption over $\Omega$ is a uniform upper bound over the distance from $x\in M$ to $\Omega$, that $\operatorname{dist}(x, \Omega)\le C$ uniformly for all $x\in M$. This can also be phrased as $\Omega$ being $\epsilon$-dense in $M$ for some large $\epsilon$. The physical heuristics behind such condition is that it is less likely for a quantum particle to tunnel into a forbidden region far away. Thus, a uniform upper bound over the distance between $\Omega$ and $M$ guarantees that particles in $M$ only need to tunnel some distance up to $C$ to get into $\Omega$. With this assumption, I proved in \cite{wan20}, the key global Carleman estimate on conic manifolds. 

\begin{theorem}[Carleman estimates on conic manifolds, \cite{wan20}]\label{t2}
Let $M$ be a conic manifold, and $\Omega$ be open set such that $\operatorname{dist}(x, \Omega)\le C$ uniformly for all $x\in M$. Let $V(x;h)\in C^\infty_b(M\times[0, 1])$. Then 
\begin{equation}\label{e2}
\|u\|_{L^2(M)}\le e^{C/h}\left(\|(h^2\Delta-V(x; h))u\|_{L^2(M)}+\|u\|_{L^2(\Omega)}\right)
\end{equation}
uniformly for all $u\in H^{2}(M)$ and $h>0$ small. 
\end{theorem}

Here we explain what a conic manifold is. Let $(\partial M, g_\partial)$ be a compact manifold of dimension $d-1$, without boundary. A conic manifold is a manifold $(M,g)$ that looks like a cone modulo a compact piece, that is there is $M_0\subset M$ compact, and $M\setminus M_0$ is $(1,\infty)_r\times \partial M$ equipped with metric
\begin{equation}
g=dr^2+r^{-2}g_\partial. 
\end{equation}
The easiest example is $\mathbb{R}^d$. Indeed, $\mathbb{R}^d\setminus \overline{B(0,1)}=(1,\infty)_r\times \mathbb{S}^{d-1}$ equipped with metric $dr^2+r^{-2}g_{\mathbb{S}^{d-1}}$. Another example are 2-dimensional cones $\{(x,y,z)\in\mathbb{R}^3: z=\alpha\sqrt{x^2+y^2}\}$ for some $\alpha\in (0, \infty)$, with its conic regularity at $0$ smoothed out. 

In this note, we will go over the proof of Theorem \ref{t2}. As an application, we will prove a logarithmic decay for Klein--Gordon waves on conic manifolds. 

\section{Local Carleman estimates}
The estimate \eqref{e2} is known to hold on $\mathbb{R}^d$ in Burq--Joly \cite{bj16}. We discuss how to obtain a estimate of type \eqref{e2} in this section, and why it is difficult to generalise it to conic manifolds. 

The standard strategy to prove a Carleman estimate on compact manifolds begins with a \emph{weight function}, $\psi\in C^\infty_b(M)$, whose properties we will specify later. We construct
\begin{equation}
\phi(x)=e^{\lambda \psi},
\end{equation}
for some $\lambda>0$ large determined later. We can assume $V(x;h)=V(x)$ only depends on $x$, since the lower order terms of the potential will not affect the estimates. Set $P_h=h^2\Delta-V(x)$, whose principal symbol is $\abs{\xi}^2-V(x)$. We then conjugate $P_h$ with exponential weights
\begin{equation}
P_{\phi}=e^{\phi/h}P_h e^{-\phi/h}=h^2\Delta-\abs{\nabla\phi}^2-V(x)+2(\nabla\phi)\cdot h\nabla-h\Delta\phi,
\end{equation}
with principal symbol
\begin{equation}
p_\phi=\abs{\xi-V(\nabla\phi)\cdot h\nabla+i\nabla\phi}_g^2-V(x)=\abs{\xi}^2-V(x)-\abs{\nabla \phi}^2+2i(\nabla\phi)\cdot \xi.
\end{equation}
We also compute the subprincipal symbol, the Poisson bracket of $\operatorname{Re}p_\phi$ and $\operatorname{Im}p_\phi$:
\begin{equation}
p_{\operatorname{sub}}=\{\operatorname{Re}p_\phi, \operatorname{Im}p_\phi\}=2(\nabla V+\nabla \abs{\nabla \phi}^2)\cdot \nabla \phi+4\xi\cdot \nabla^2\phi \cdot \xi. 
\end{equation}
The goal is to establish that, $p_\phi=0\Rightarrow p_{\operatorname{sub}}\neq 0$ where $\nabla\phi\neq 0$. Indeed, let $\chi\in C^\infty(M)$ be a cutoff function such that $\abs{\nabla\phi}\ge \rho$ for a small constant $\rho>0$ on $\supp\chi$. Then on the neighbourhood of the characteristic variety of $P_\phi$, 
\begin{equation}
\Sigma=\left\{\frac{1}{4}(\abs{\nabla \phi}^2+V(x))\le \abs{\xi}^2\le 4(\abs{\nabla \phi}^2+V(x))\right\},
\end{equation}
intersecting $\supp\chi$, we have
\begin{equation}
\frac{1}{4}(\lambda^2\rho^2e^{2\lambda \psi}+V(x))\le \abs{\xi}^2\le 4(C\lambda^2e^{2\lambda \psi}+V(x))
\end{equation}
and furthermore via direct computation that after choosing a large $\lambda$, we can have
\begin{equation}
p_{\operatorname{sub}}\ge 2\rho^4\lambda^4 e^{3\lambda\psi}+\bigo(\lambda^3 e^{3\lambda\psi})\ge \epsilon>0
\end{equation}
uniformly. The fact that the principal and subprincipal symbols cannot both vanish over the $\supp\chi$ means $P_\phi$ is subelliptic: via the Gårding inequality we have
\begin{equation}
\|P_\phi \chi u\|_{L^2}\ge Ch^{\frac{1}{2}}\|\chi u\|_{H_h^{\frac{3}{2}}}. 
\end{equation}
Note $P_\phi=e^{\phi/h}P_he^{-\phi/h}$ to see for any $\chi$ with $\abs{\nabla \phi}\ge \rho$ on the $\supp\chi$, we have the key local estimates
\begin{equation}\label{e4}
\|e^{\phi/h}\chi u\|_{L^2}\le Ch^{-\frac{1}{2}}\|e^{\phi/h}P_h\chi u\|_{L^2}.
\end{equation}
This reduces, via taking minimum and maximum of $e^{\phi/h}$, to
\begin{equation}
\|\chi u\|_{L^2}\le C e^{C/h}(\|P_h u\|_{L^2}+\|[P_h, \chi] u\|_{L^2})
\end{equation}
Now we assume that $\psi$ is a function such that $\abs{\nabla \phi}\ge\rho$ on $M\setminus \Omega$. Then we can choose $(1-\chi)$ supported in $M\setminus \Omega$, and have
\begin{gather}
\|(1-\chi)u\|_{L^2}\le C\|\mathbbm{1}_{\Omega} u\|_{L^2}\\
\|[P_h, \chi] u\|_{L^2}\le C\|P_h u\|_{L^2}+C\|\mathbbm{1}_{\Omega} u\|_{L^2}. 
\end{gather}
We eventually have
\begin{equation}
\|u\|_{L^2}\le e^{C/h}(\|P_h u\|_{L^2}+\|\mathbbm{1}_{\Omega}u\|_{L^2})
\end{equation}
as desired. 

Therefore, we can always have the desired Carleman estimates \eqref{e2}, as long as we can find a weight function $\psi\in C^\infty_b(M)$ such that $\abs{\nabla \phi}\ge \rho$ on $M\setminus \Omega$. We here discuss the strategy to find such a function. Firstly, on a compact manifold $M$, we can pick a Morse function $\psi$ (a function with finitely many critical points), and then diffeomorphically pull those critical points back into $\Omega$ one by one. Such function $\psi$ after finitely many diffeomorphism will be an ideal function. We thus know the Carleman estimate \eqref{e2} indeed holds for any compact manifold $M$ and arbitrary open set $\Omega$. 

We now consider the case when $M=\mathbb{R}^d$ is noncompact, $\Omega$ is $\epsilon$-dense considered in \cite{bj16}. It is natural to again try to find a weight function whose critical points are also $\epsilon$-dense in $\mathbb{R}^d$. On $\mathbb{R}^d$,
\begin{equation}
\psi(x)=\cos(2\pi x)=\Pi_{k=1}^d \cos(2\pi x_k)
\end{equation}
has critical points on lattices $\mathbb{Z}^d\subset\mathbb{R}^d$. This was the way how a Carleman estimate was proved in \cite{bj16}. 

In our case of conic manifolds, there does not exist a good lattice structure like $\mathbb{Z}^d$, such that the lattice points are $\epsilon$-dense in $M$ and not too close to each other. Even with such lattice structure, there might be no functions whose critical points are exactly at those lattices. Such lack of structure encourages us to find a new approach to the Carleman estimates. 

\section{Carleman estimates with multiple weights}
In \cite{bur98}, Burq worked on a Carleman estimate for systems where topological constraints prevent any weight function on compact $M$ to have all critical points inside some specific open set $\Omega$. In order to address the failure to construct one suitable weight function for the Carleman estimate, Burq introduced two weight functions $\psi_+, \psi_-$ that are compatible with each other in the following sense:
\begin{equation}\label{e3}
\nabla \psi_\pm(x)=0, x\notin \Omega \Rightarrow  \nabla\psi_{\mp}(x)\neq 0, \psi_{\mp}(x)> \psi_{\pm}(x).
\end{equation}
It is interpreted that the weight functions may have critical points outside $\Omega$, but at each critical point outside $\Omega$, the other weight function must not have a critical point and must dominate the weight function whose gradient vanishes. With two compatible weight functions, Burq showed in \cite{bur98} that one indeed can recover \eqref{e2}. 

I was motivated by the power of the two-weight Carleman argument to overcome topological constraints: can one similarly devise of a method of divide-and-conquer to assembly local Carleman estimates on finitely many parts of $M$ into a global estimate \eqref{e2}? In order to do so, we need a stronger Carleman estimates that accommodate more than two compatible weights, that is, a \emph{Carleman estimate with multiple weights}. 

\begin{definition}[Compatibility]
We say a finite family of weight functions $\{\psi_k\}_{k=1}^n$ in $C_b^{\infty}(M)$ is compatible with the control from $\Omega$, if there are constants $\rho, \tau>0$ such that, at each point $x\in M\setminus \Omega$, if $\abs{\nabla \psi_k}<\rho$ for some $k$, then there exists some $l$ such that
\begin{equation}
\abs{\nabla \psi_l}\ge \rho, \ \psi_l(x)\ge \psi_k(x)+\tau. 
\end{equation}
\end{definition}

The compatibility condition can be see as a natural generalisation of \eqref{e3}: if one of the weight $\psi_k$ has small gradient outside $\Omega$, then one of the other weights $\psi_l$ must have large gradient and dominate $\psi_k$ thereat. With the compatibility condition, we can prove the key Carleman estimate with multiple weights:

\begin{theorem}[Global Carleman estimates with multiple weights, \cite{wan20}]\label{t1}
Let $M$ be a manifold of bounded geometry and $V(x;h)\in C^\infty_b(M\times[0,1])$. Assume that there are weights $\{\psi_k\}_{k=1}^n$ compatible with the control from $\Omega$. Then there is $C>0$ such that
\begin{equation}
\|u\|_{L^2(M)}\le e^{C/h}\left(\|(h^2\Delta-V(x; h))u\|_{L^2(M)}+\|u\|_{L^2(\Omega)}\right)
\end{equation}
uniformly for all $u\in H^{2}(M)$ and $h>0$ small. 
\end{theorem}
Here, \emph{manifolds of bounded geometry} are defined by Shubin in \cite{shu92}: they are manifolds on which all derivatives of the Riemannian metric tensor are uniformly bounded, and the injectivity radius is bounded from below. Conic manifolds are manifolds of bounded geometry. The new tool, Carleman estimates with multiple weights, allow us to cut the manifold up into finitely many pieces, and glue local Carleman estimates on each piece up to form a global estimate. Such gluing argument is not trivial in semiclassical analysis: indeed, the $e^{C/h}$ factor will make any $h^N$-order error exponentially large. 

The rest of this section is devoted to the proof of Theorem \ref{t1}. Firstly consider cutoffs $\chi_k\in C^\infty_c(M)$ such that $\chi_k\equiv 1$ when $\abs{\nabla \psi_k}\ge \rho$: heuristically, $\chi_k$ localises where a local Carleman estimate \eqref{e4} with weight $\psi_k$ can control. We observe that $\chi_k$ can be chosen that $\sum_k \chi_k\ge 1$ on $M\setminus\Omega$ due to the compatibility condition. 

We resume the analysis from the local estimates \eqref{e4}: for each $\chi_k$, let $\phi_k=e^{\lambda\psi_k}$ and we have
\begin{equation}
\|e^{\phi_k/h}\chi_k u\|_{L^2}\le Ch^{-\frac{1}{2}}\|e^{\phi_k/h}P_h\chi_k u\|_{L^2}.
\end{equation}
Using similar estimates on the commutators as before, we can show that
\begin{equation}
\|e^{\phi_k/h}\chi_k u\|_{L^2}\le Ch^{-\frac{1}{2}}\|e^{\phi_k/h}P_h u\|_{L^2}+Ch^{-\frac{1}{2}}\|e^{\phi_k/h}\tilde\chi_k u\|_{L^2}+Ch^{-\frac{1}{2}}\|e^{\phi_k/h}\mathbbm{1}_\Omega u\|_{L^2},
\end{equation}
where $\supp\tilde\chi_k\subset \{\nabla \psi_k<\rho\}\setminus \Omega$. We then sum up the local estimates to have a crude global estimate:
\begin{equation}\label{e5}
\sum_k\|e^{\phi_k/h}\chi_k u\|\le C\sum_k\left(h^{-\frac{1}{2}}\|e^{\phi_k/h}P_h u\|+h^{-\frac{1}{2}}\|e^{\phi_k/h}\tilde\chi_k u\|+h^{-\frac{1}{2}}\|e^{\phi_k/h}\mathbbm{1}_\Omega u\|\right).
\end{equation}

We now improve the estimate \eqref{e5} from the left and right. The compatibility condition guarantees the inequality 
\begin{equation}
 \sum_k e^{\phi_k/h}\chi_k\ge \frac{1}{2n}\sum_k e^{\phi_k/h}
 \end{equation} 
on $M\setminus \Omega$. Thus from the left of \eqref{e5}, we have
\begin{equation}
\frac{1}{2n}\|(\sum_ke^{\phi_k/h})(1-\mathbbm{1}_\Omega) u\|\le \sum_k\|e^{\phi_k/h}\chi_k u\|.
\end{equation}
On another hand, the compatibility condition gives the inequality 
\begin{equation}
e^{\phi_{k'}/h}\tilde\chi_{k'}\le e^{-\epsilon/h}\sum_k e^{\phi_k/h}, \ \epsilon=(e^{\lambda\tau}-1)e^{\lambda \min_{k}(\inf_M\psi_k)}>0,
\end{equation}
and an bound for the second term on the right of \eqref{e5}
\begin{equation}
\sum_{k}\|e^{\phi_{k}/h}\tilde\chi_{k}u\|\le nh^{-\frac{1}{2}}e^{-\epsilon/h}\|(\sum_k e^{\phi_k/h})u\|.
\end{equation}
Thus \eqref{e5} is improved to 
\begin{multline}
\|(\sum_ke^{\phi_k/h})(1-\mathbbm{1}_\Omega) u\|\\
\le C\sum_k\left(h^{-\frac{1}{2}}\|e^{\phi_k/h}P_h u\|+h^{-\frac{1}{2}}e^{-\epsilon/h}\|(\sum_k e^{\phi_k/h})u\|+h^{-\frac{1}{2}}\|e^{\phi_k/h}\mathbbm{1}_\Omega u\|\right).
\end{multline}
Add $\|(\sum_ke^{\phi_k/h})\mathbbm{1}_\Omega u\|$ to the left and absorb the second term on the right to see
\begin{equation}
|(\sum_ke^{\phi_k/h})u\|\le C\sum_k\left(h^{-\frac{1}{2}}\|e^{\phi_k/h}P_h u\|+h^{-\frac{1}{2}}\|e^{\phi_k/h}\mathbbm{1}_\Omega u\|\right).
\end{equation}
Now take the minimum and maximum over both sides to see
\begin{equation}
\|u\|_{L^2(M)}\le e^{C/h}\left(\|P_h u\|_{L^2(M)}+\|u\|_{L^2(\Omega)}\right).
\end{equation}
Here we conclude the proof of Theorem \ref{t1}.

\section{Weight functions construction}
What remains to prove Theorem \ref{t2} is a family of compatible weight functions on conic manifolds. In this section we briefly explain how to construct them. 

\begin{figure}[h]
\centering{
\includegraphics{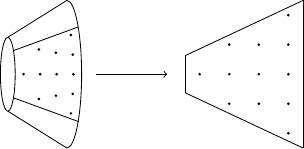}}
\caption{\label{f1}Cutting up conic ends into sectors, where the dots stand for the critical points of weight functions.}
\end{figure}

We begin with the conic part $M\setminus M_0$. We cut them into finitely many curved sectors invariant along $r>C$. We then quasi-isometrically embed each curved sector into flat sectors in $\mathbb{R}^d$. We then pull back the $\cos(2\pi x)$ function localised in the flat sectors back onto weight functions $\psi_k$ on curved sectors. We move individual critical points into $\Omega$ one by one via diffeomorphism. Carefully trimming the edges of $\psi_k$, we can make sure that the weights are compatible with each other in the overlapped regions between different curved sectors. See Figure \ref{f1} for illustration. 

For the central compact part $M_0$, we can pick a Morse function and repeat the same strategy on a compact manifold. What is extra is that we need to trim the resulted weight function near the boundary of $M_0$ such that it is compatible with other weight functions on the conic part. 

It is technical to check, but such construction is possible: see further details in \cite{wan20}. Thus we have proved Theorem \ref{t2} by establishing a family of compatible weight functions for Theorem \ref{t1}. 

\section{Application: logarithmic decay}
As an application of Theorem \ref{t2}, we prove a logarithmic decay result for damped Klein--Gordon equations. Consider the damped Klein--Gordon equations on conic manifolds, 
\begin{gather}
(\partial_t^2+W(x)\partial_t+\Delta+1)u(t,x)=0, \ (t,x)\in \mathbb{R}_{t\ge 0}\times M,\\
u(0, x)=u_0(x)\in H^2(M), \ \partial_t u(0, x)=u_1(x)\in H^1(M).
\end{gather}
Here $W(x)\in C^\infty_b(M)$ is the damping function, measuring the dissipation of waves on $M$. The energy of the Klein--Gordon waves is defined by 
\begin{equation}
E(u,t)=\frac{1}{2}(\|\nabla_x u\|_{L^2(M)}^2+\|\partial_t u\|_{L^2(M)}^2),
\end{equation}
and it decays due to the damping $W(x)\ge 0$:
\begin{equation}
\partial_t E(u,t)=-\|\sqrt{W(x)}\partial_t u\|_{L^2(M)}^2\le 0. 
\end{equation}
We have the following corollary of Theorem \ref{t2}, our Carleman estimate on conic manifolds:
\begin{corollary}[Logarithmic decay, \cite{wan20}]\label{t3}
Assume $\{W\ge \delta\}$ is $\epsilon$-dense in $M$ conic, for some $\delta,\epsilon>0$. Then there exists $C>0$ such that
\begin{equation}
E(u,t)\le \frac{C}{\log(2+t)}(\|u_0\|_{H^2(M)}^2+\|u_1\|_{H^1(M)}^2)
\end{equation}
uniformly for any solution $u$ to the Klein--Gordon equations with initial data $(u_0, u_1)$. 
\end{corollary}
The proof of Corollary \ref{t3} reduces to proving an uniform estimate of the form
\begin{equation}
\|(h^2\Delta-ihW(x)-1)^{-1}\|_{\mathcal{L}(L^2)}\le e^{C/h}
\end{equation}
for all $h$ small, using the semigroup stability results in \cite{bur98}. Let 
\begin{equation}
(h^2\Delta-ihW(x)-1)u=f.
\end{equation}
Since $\{W\ge \delta\}$ is $\epsilon$-dense, apply Theorem \ref{t2} that
\begin{multline}\label{e6}
\|u\|\le e^{C/h}(\|(h^2\Delta-1)u\|+\|\mathbbm{1}_{\{W\ge \delta\}}u\|)\le e^{C/h}\|(h^2\Delta-1)u\|+e^{C/h}\|\sqrt{W}u\|\\
\le e^{C/h}\|f\|+Ce^{C/h}\|\sqrt{W}u\|,
\end{multline}
where we used $\|Wu\|\le C\|\sqrt{W}u\|$. Pair $f$ with $u$ to see
\begin{equation}
\langle f,u\rangle=\|h\nabla u\|^2-\|u\|^2-ih\|\sqrt{W}u\|^2,
\end{equation}
whose imaginary part implies
\begin{equation}
\|\sqrt{W}u\|\le 4\epsilon^{-1}h^{-1}e^{C/h}\|f\|+\epsilon e^{-C/h}\|u\|.
\end{equation}
Revisit \eqref{e6} and absorb the $\epsilon$-term to see
\begin{equation}
\|u\|\le e^{C'/h}\|f\|,
\end{equation}
giving the desired exponential bound for the resolvent. The proof of Corollary \ref{t3} is now concluded.

\bibliographystyle{alpha}
\bibliography{Robib}

\end{document}